\documentclass[11pt]{article}

\usepackage{graphicx}

\newcommand{\Du}{{D\" urer}}

\begin{document}

\title{Party Game for a 500th Anniversary}

\author{Fumiko Futamura,$^{1,}$\thanks{futamurf@southwestern.edu}  \ 
Marc Frantz,$^{1,}$\thanks{marcfrantz10@gmail.com}\ \ 
and 
Annalisa Crannell$^{2,}$\thanks{annalisa.crannell@fandm.edu}\\[5pt]
$^1${\itshape Southwestern University, Georgetown, TX 78626, USA}\\[5pt]
$^2${\itshape Franklin \& Marshall College, Lancaster, PA 17604, USA}}

\date{\today}

\maketitle

\begin{abstract}
On the 500th anniversary of Albrecht \Du's copperplate engraving {\it Melencolia~I}, we invite readers to join in  a time-honored ``party game"  that has  attracted art historians and scientists for many years: guessing the nature and meaning of the composition's enigmatic stone polyhedron. Our main purpose is to demonstrate the usefulness of the cross ratio in the analysis of works in perspective. We show how the cross ratio works as  a projectively invariant ``shape parameter" of the polyhedron, and how it can be used in analyzing various theories of this figure. \newline
\null\newline
\null\newline
\end{abstract}

\section{Introduction}

\begin{figure}[here!]
\centering
\raisebox{0.9in}{\includegraphics[width=1.5in]{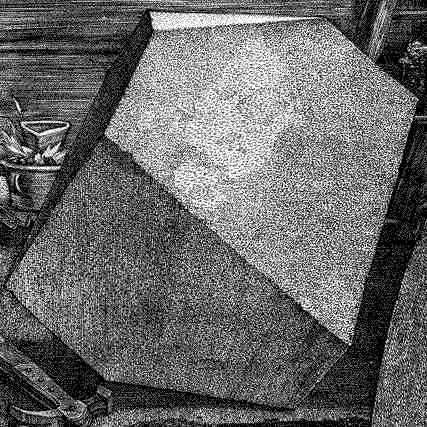}}
\qquad
\includegraphics[scale=1]{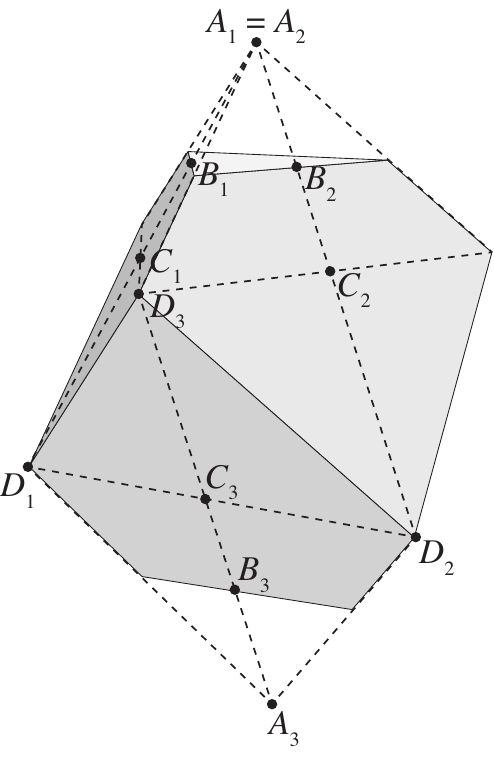} 
\caption{Detail of {\it Melencolia~I\/} (left) and a blowup of the polyhedron.}
\label{setup}
\end{figure}

This year marks the 500th anniversary of the copperplate engraving {\it Melencolia~I\/}  by the great artist-mathematician Albrecht \Du\ (1471--1528). A good way for math lovers to celebrate this anniversary is to play a time-honored ``party game"  that has  attracted art historians and scientists for many years: guessing the nature and meaning of the composition's enigmatic stone polyhedron, illustrated in Figure~\ref{setup}\footnote{A high-resolution image of the complete engraving, convenient for the purpose of analysis, can be found at http://upload.wikimedia.org/wikipedia/commons/1/18/D\" urer\_Melancholia\_I.jpg.}. Published attempts at playing the game date back at least a century \cite{Weitzel}. In fact, it's a game that may go on indefinitely, because (1) we have no writings by \Du\ definitively saying what the polyhedron is, (2) our ability to measure and analyze the polyhedron is subject to at least some small errors, and (3) the image itself has anomalous features, some of which we discuss in the appendix. 

Nevertheless the polyhedron, sometimes known as ``\Du's solid", seems to be accurately drawn for the most part, and authors such as MacGillavry \cite{MacGillavry} have been able to get good fits of their models by allowing for a certain amount of human error in both drawing and measuring. In this article we discuss a possible model of the solid and compare it to the results of MacGillavry and others. We will try to make the case that one of the most convenient and effective ways to score the game is to use the cross ratio. We show how the cross ratio works as a projectively invariant ``shape parameter" of the polyhedron, and how it can be used in analyzing various theories of this figure. 

\section{Playing the game}

A commonly accepted rule of the game is the assumption that the polyhedron is formed by starting with a {\it trigonal trapezohedron} or {\it trigonal deltohedron}---a three-dimensional figure whose faces are six congruent rhombi (a cube is such a figure).  Then, with the longest diagonal held vertically, congruent tetrahedrons are cropped off the top and bottom by horizontal planes, leaving congruent equilateral triangles as the top and bottom  faces.  The other six faces---the truncated rhombi---are thus congruent pentagons. The shape of the solid is completely determined by the shape of one of the pentagonal faces, hence many authors describe their model of the solid by describing the shape of one of these pentagons.

In addition to rules, the game should have some method for scoring various attempts at it, at least in an informal way. First of all, it makes no sense to play the game if we assume that \Du\ was such a poor draftsman that we can simply ignore the engraving and conjecture any shape we want. We therefore require that a proposed model of the solid should have a shape that is  close to that derived by perspective analysis of the engraving.  In addition to respecting \Du\ as an artist, we must acknowledge his reputation as a mathematician. In the words of the respected mathematician Morris Kline, ``\ldots of all the Renaissance artists, the best mathematician was the German Albrecht \Du \ldots " \cite[p.~233]{Kline}. Among other things, \Du\ owned a copy of Euclid and wrote well-respected treatises on proportion, including---significantly for us---the golden ratio. Thus we assume that \Du\ incorporated some kind of interesting mathematical relationships into the design of the solid, rather than just sketching until he found something that looked pleasing to him. A proposed model of the solid should include the description of such relationships.

\section{A useful parameter: the cross ratio}
\begin{figure}[here!]
\centering
\includegraphics{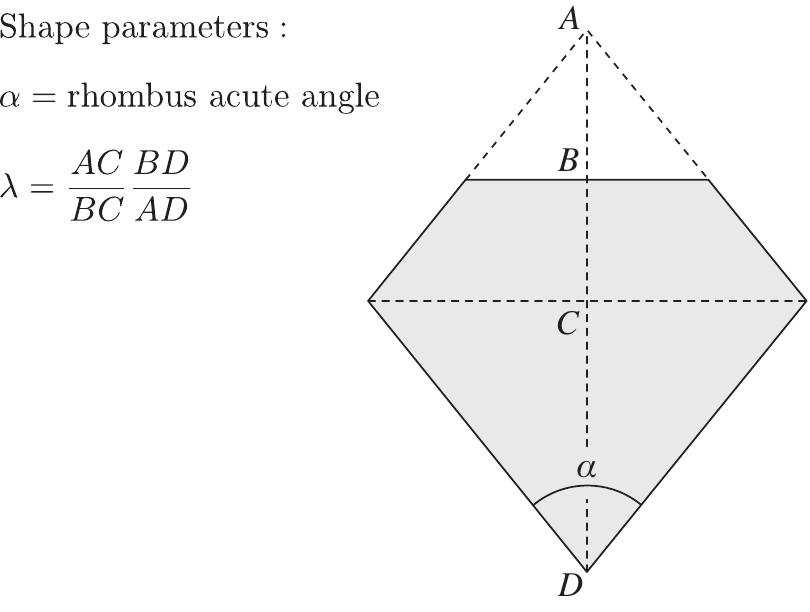}
\caption{The parameters $\alpha$ and $\lambda$ determine the shape of the truncated rhombic face and hence the shape of the solid.}
\label{generic-face}
\end{figure}

Under the assumptions of the game, the shape of the solid is completely determined by two numbers, which we refer to as ``shape parameters". As we mentioned earlier, the shape of the solid is determined by the shape of one of the six congruent pentagonal faces. The shape of a pentagonal face (see Figure~\ref{generic-face}) can be specified by the acute angle $\alpha$ of the rhombus and another parameter $\lambda$, which determines the level at which the rhombus is truncated. One obvious candidate for $\lambda$ would be the ratio of distances $BC/AC$ in Figure~\ref{generic-face}; that is, the fraction of $AC$ left after truncation. However, this parameter has one disadvantage; it cannot always be measured directly in a perspective drawing of the face, because perspective often distorts the ratios of lengths.

Fortunately, given four collinear points $A$, $B$, $C$, $D$ as in Figure~\ref{generic-face}, there is a quantity whose value is not changed by perspective, namely the {\it cross ratio}\footnote{The cross ratio is sometimes denoted by $(AB,CD)$, among other things.} $\times(ABCD)$ defined by
\[
\times(ABCD)=\frac{AC}{BC}\frac{BD}{AD}.
\]
For a truncated rhombic face like the one in  Figure~\ref{generic-face}, there is a convenient relationship between the cross ratio and the  ratio $BC/AC$. Letting $\lambda=\times(ABCD)$, we have
\begin{equation}\label{lambda}
\lambda=\times(ABCD)=\frac{AC}{BC}\frac{BD}{AD}=\frac{AC}{BC}\frac{BC+AC}{2AC}=\frac{1}{2}\left(1+\frac{1}{BC/AC}\right) .
\end{equation}
Solving (\ref{lambda}) for $BC/AC$ gives
\begin{equation}\label{BC/AC}
\frac{BC}{AC}=\frac{1}{2\lambda-1} .
\end{equation}

From (\ref{lambda}) and (\ref{BC/AC}) we see that the truncation of the rhombus determines the value of the cross ratio $\lambda$, and the value of $\lambda$ determines the truncation of the rhombus. Since $\lambda$ is projectively invariant---that is, unchanged by the distortion of perspective---we choose it as the second shape parameter.

\begin{figure}[here!]
\centering
\includegraphics{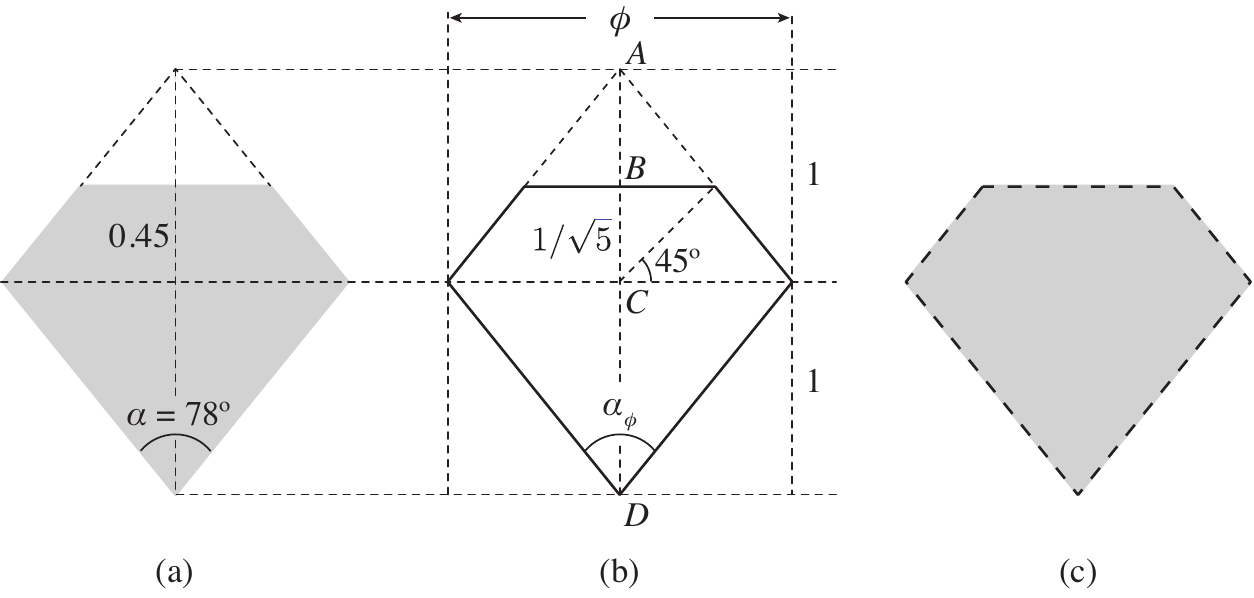}
\caption{In (a) the $78^\circ$ version of MacGillavry's pentagonal face in solid gray. In (b) the ``golden" pentagon in black outline, with the rhombus inscribed in a pair of golden rectangles. In (c) MacGillavry's pentagon in solid gray, and the golden pentagon in dashed black outline. When drawn at this size, they are practically identical.}
\label{m-face}
\end{figure}

Perhaps the most thorough and frequently cited perspective analysis of the solid is the one done by MacGillavry \cite{MacGillavry}.  MacGillavry estimated the shape parameters of the pentagonal faces by computing them several different ways, applying the rules of perspective to different parts of the engraving and averaging the results. MacGillavry estimated $\alpha=79^\circ\pm 1^\circ$ and $BC/AC\approx 0.45$.  Figure~\ref{m-face}(a) shows such a pentagon in solid gray, with MacGillavry's minimum value of $\alpha=78^\circ$. From this value of $BC/AC$ and equation (\ref{lambda}) we see that MacGillavry's work implies
\[
\lambda=\frac{1}{2}\left( 1+\frac{1}{.45}\right)\approx 1.61.
\]

Our own estimate of $\lambda$ is close to that of MacGillavry. In Figure~\ref{setup} we have extended sides of the three visible pentagonal faces to locate the truncated vertices $A_1\ (=A_2)$ and $A_3$ of the rhombi. For $n=1,2,3$, the lines $A_nD_n$ are the centerlines of the rhombi, the points $C_n$ are the centers, and the points $B_n$ are the truncation points of the centerlines. If our assumptions about the polyhedron are correct, and if the drawing and our measurements were perfect, then the cross ratios $\times(A_nB_nC_nD_n)$ would be exactly the same for $n=1,2,3$. In an imperfect world, we of course expect some variation. One of the authors computed these cross ratios by importing a digital image of {\it Melencolia~I\/} into the program GeoGebra, and found 
\begin{equation}\label{crosses}
\times(A_1B_1C_1D_1)=1.63,\quad \times(A_2B_2C_2D_2)=1.64,\quad
\times(A_3B_3C_3D_3)=1.59.
\end{equation}
To allow for imperfections as MacGillavry did, we average these three values to get $1.62$. 

The reader may have noticed that these estimates of $\lambda\approx 1.61$  and $\lambda\approx 1.62$  are suspiciously close to a famous number: the golden ratio $\phi$, where
\[
 \phi=\frac{\sqrt{5}+1}{2}\approx 1.618.
 \]
Following this clue, our model of the solid is based on the golden ratio, with parameters we denote by $\alpha_\phi$ and $\lambda_\phi$. Figure~\ref{m-face}(b) shows our proposed model of a face in black outline. The rhombus is inscribed in a pair of golden rectangles, each with height $1$ and width $\phi$. The truncation line is easily found by drawing a $45^\circ$ line from the center of the rhombus until it meets an edge of the rhombus as shown. It is straightforward to show that $BC/AC=BC=1/\sqrt{5}$, resulting in a cross ratio, via (\ref{lambda}), of 
\[
\lambda_\phi=\frac{1}{2}\left( 1+\frac{1}{1/\sqrt{5}}\right)=\frac{\sqrt{5}+1}{2}=\phi .
\]

Thus the cross ratio is the golden ratio. The acute angle $\alpha_\phi$ is given by  
\[
\alpha_\phi =2\arctan(\phi/2)\approx 77.95^\circ.
\]
The values of  $\lambda_\phi$ and $\alpha_\phi$ are very close to MacGillavry's values of $\lambda=1.61$ and $\alpha=78^\circ$. In fact, as shown in Figure~\ref{m-face}(c),  the two pentagons are nearly indistinguishable (MacGillavry's is solid gray and ours is the dashed black outline). For the above reasons  we (somewhat playfully) refer to the pentagon in (b) as a {\it golden pentagon}.

To analyze how close the two pentagons really are, observe that the half-width of the pentagon in Figure~\ref{m-face}(a), divided by the half-width of the pentagon in (b) is
$
\tan(78^\circ/2)\div(\phi/2)\approx 1.00095
$.
Thus if the half-widths of the pentagons are about an inch (they're actually a little less), then the difference in the half-widths is less than a thousandth of an inch. The ratios of the heights of each pentagon above the horizontal centerline  is
$
0.45\div(1/\sqrt{5})\approx 1.0062
$.
Thus if these heights are about a half inch (again they're a little less), then the difference between them is about three thousandths of an inch. We conclude that from an intuitive, visual standpoint, MacGillavry's minimum-angle solution of $\lambda=1.61$ and $\alpha=78^\circ$ is essentially the same as ours, minus the ``golden" formulation of the shape parameters: $\lambda=\phi$ and $\alpha=2\arctan(\phi/2)$.

\section{Testing theories with the cross ratio}

Over the years many theories have been advanced as to the shape of \Du's solid, or equivalently, its pentagonal faces. With so many possibilities, any particular one of them (including this one) probably stands only a small chance of being right. The advantage of this abundance of theories is that it can be as much fun testing old theories as it is cooking up new ones. Following our methods of the previous section, we show how the cross ratio can be applied to this variation of the game. 

\subsection{Lynch's model}
\begin{figure}[here!]
\centering
\includegraphics{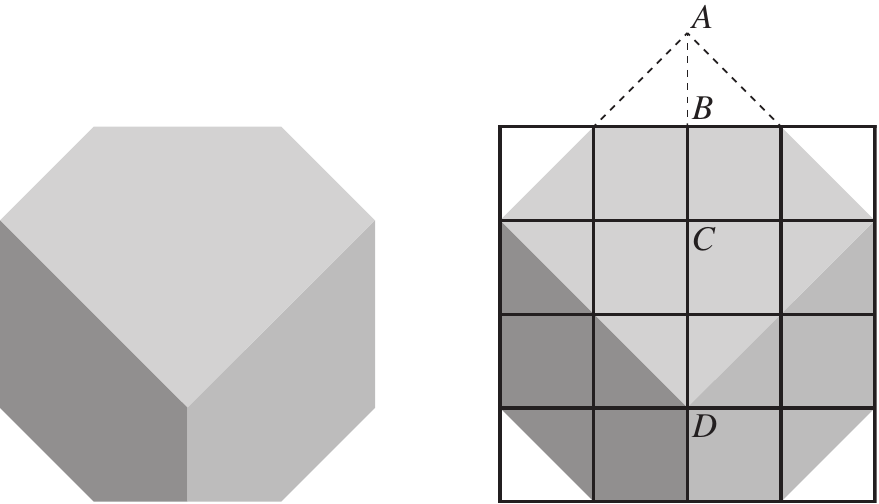}
\caption{Diagram after a figure by Lynch \cite[p.~228]{Lynch}.}
\label{4-grid}
\end{figure}

Figure~\ref{4-grid} depicts a model conjectured by Lynch \cite{Lynch}, who guessed that vertices of the solid project orthogonally onto the lattice points of a $4\times 4$ grid, thus linking it with a $4\times 4$ magic square that appears in {\it Melencolia~I}.  However, if this were true, we would have (assuming for simplicity that the squares are 1 unit on a side),
\[
\lambda=\times(ABCD)=\frac{AC}{BC}\frac{CD}{AD}=\frac{2}{1}\frac{3}{4}=\frac{3}{2},
\]
which differs significantly (about $7\%$) from our measured value of $1.62$ and MacGillavry's estimate of $1.61$.   We conclude that Lynch's idea is not in good agreement with \Du's perspective rendering of the solid.

\subsection{Schreiber's model}
\begin{figure}[here!]
\centerline{\includegraphics{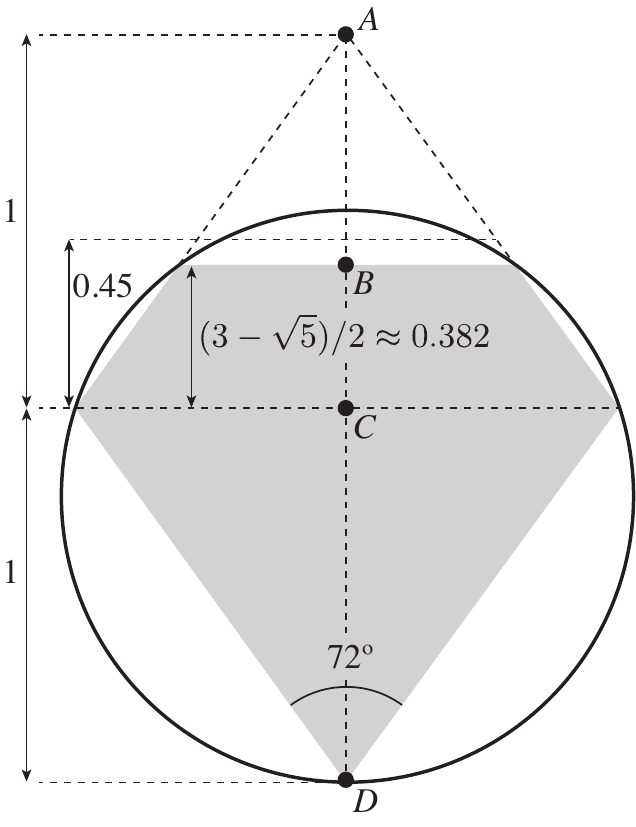}}
\caption{The $72^\circ$, circle-inscribed face conjectured by Schreiber \cite{Schreiber}.}
\label{72}
\end{figure}

Figure~\ref{72} depicts a model by Schreiber having an acute rhombus angle of $72^\circ$ and a truncation chosen so that the face is inscribed in a circle. But a straightforward computation shows that this requires a ratio of $BC/AC=(3-\sqrt{5})/2$. When substituted into (\ref{lambda}), this gives a cross ratio value of $\lambda\approx 1.81$, which differs from MacGillavry's estimate and ours by more than $11\%$. This can be checked approximately by using a ruler and protractor to draw the $72^\circ$ rhombus, and a circle template to choose a circle that passes through its left, right, and bottom vertices. The circle then determines the truncation. In either case the cross ratio of the resulting pentagon contradicts the evidence of the engraving itself. In addition, the $72^\circ$ angle is $8\%$ less than MacGillavry's minimum value of $78^\circ$.

\subsection{Weitzel's model}
\begin{figure}[here!]
\centerline{\includegraphics{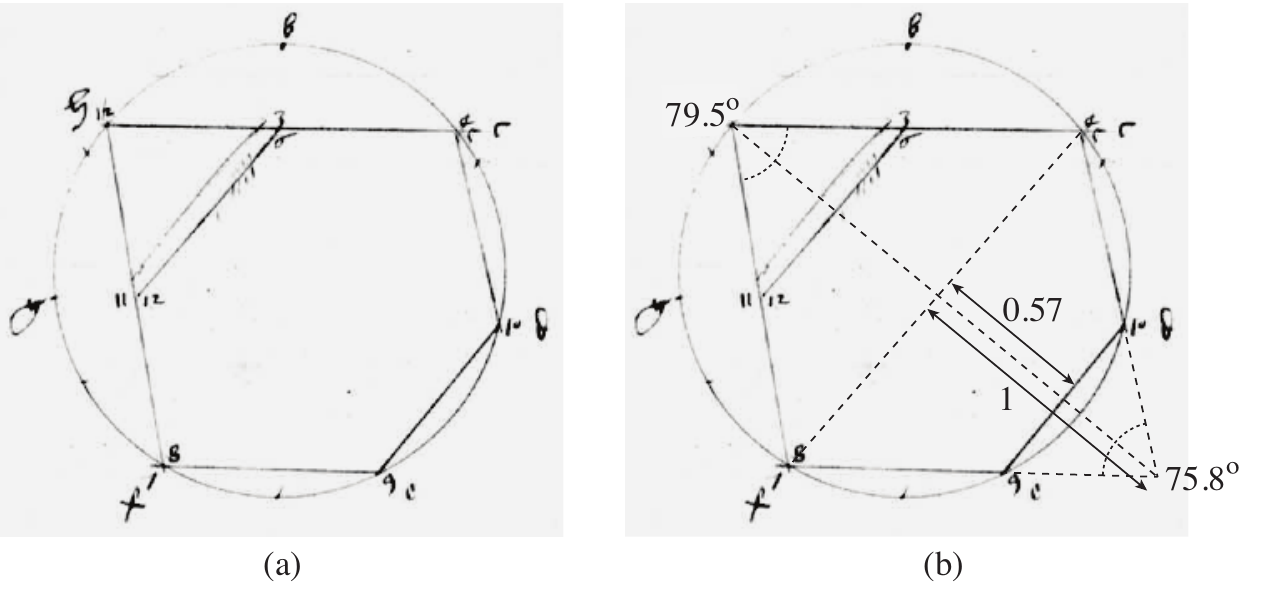}}
\caption{In (a), \Du's sketch discussed by Weitzel \cite{Weitzel}. In (b), some measurements of angles and lengths.}
\label{sketch}
\end{figure}

Figure~\ref{sketch} shows a drawing from a \Du\ sketchbook found by Weitzel \cite{Weitzel}. Weitzel conjectured that this drawing is a preliminary sketch of a face of the solid in {\it Melencolia~I\/}, and estimated the rhombus angle to be $79.5^\circ\pm 0.5^\circ$, apparently based on the angle at the upper left in Figure~\ref{sketch}(b). However, we noticed that the drawing is not perfectly symmetrical, as the other angle is about $75.8^\circ$. The average of these values is close to MacGillavry's lower bound of $78^\circ$. A feature of the sketch that is {\it not\/} consistent with MacGillavry's results is the truncation ratio of approximately $0.57$ indicated in Figure~\ref{sketch}(b), apparently intended to allow for a circumscribed circle.  This is significantly different from MacGillavry's value of $0.45$ for the faces of \Du's solid. Alternatively, letting $BC/AC=.57$ in equation (\ref{lambda}) yields a cross ratio $\lambda$ for Weitzel's model of
\[
\lambda=\frac{1}{2}\left( 1+\frac{1}{.57}\right)\approx 1.38,
\]
which is again significantly different from our measured value of $1.62$ and MacGillavry's estimate of $1.61$. This suggests that \Du\ may have used an outer rhombus much like that in Figure~\ref{sketch} for the face of his solid, but truncated it differently for some reason. Assuming that the sketch really is for {\it Melencolia~I\/}, our suggestion for the different truncation would be that the resulting face has a cross ratio equal to the golden ratio.

\section{Conclusion}
Under the usual assumptions, the shape of \Du's  solid can be specified by a pair $(\alpha,\lambda)$ of shape parameters, as illustrated in Figure~\ref{generic-face}. Of the two, the cross ratio $\lambda$ is the easiest to investigate empirically, because it is projectively invariant. The two parameters, which are independent of one another, are equally important in determining the solid, but in the literature the angle $\alpha$ seems to have attracted more attention. For example, Ishizu \cite{Ishizu} provides a table of values of $\alpha$ proposed by various authors going back to 1900, but no corresponding table of values of $\lambda$, or any other parameter that would determine the truncation of the rhombus. We urge future players of the game to give equal consideration to both parameters, and we suggest the projectively invariant cross ratio as the parameter to determine the truncation.

Interestingly, Ishizu's table shows that in the literature, the proposed values of $\alpha$ cluster around either $80^\circ$ (largely because it is consistent with the perspective of the engraving) or $72^\circ$ because of its connection with the golden ratio---a connection that might have attracted  \Du, namely
\[
\alpha=2\arccos(\phi/2).
\] 
But if that is a nice connection with the golden ratio---one that would appeal to \Du\ the mathematician---then surely the double connection $(\alpha_\phi,\lambda_\phi)$ which we propose,
\[
\alpha_\phi=2\arctan(\phi/2),\quad \lambda_\phi=\phi ,
\]
must have some appeal as well, especially since it is compatible with the perspective of \Du's rendering of the solid.

Fortunately for the anniversary celebration, neither this model nor any other will end the game of guessing the intended shape of the solid. Any proposed values of $\alpha$ and $\phi$ can agree only approximately with the engraving and, as we show in the appendix, the engraving itself has some anomalous features that make measurements of it even more uncertain. There's plenty of room for more ideas and more players. We invite the reader to play the game and invent a new and better theory!

\section*{Appendix: Anomalous features of the engraving}
\begin{figure}[here!]
\centering
\includegraphics{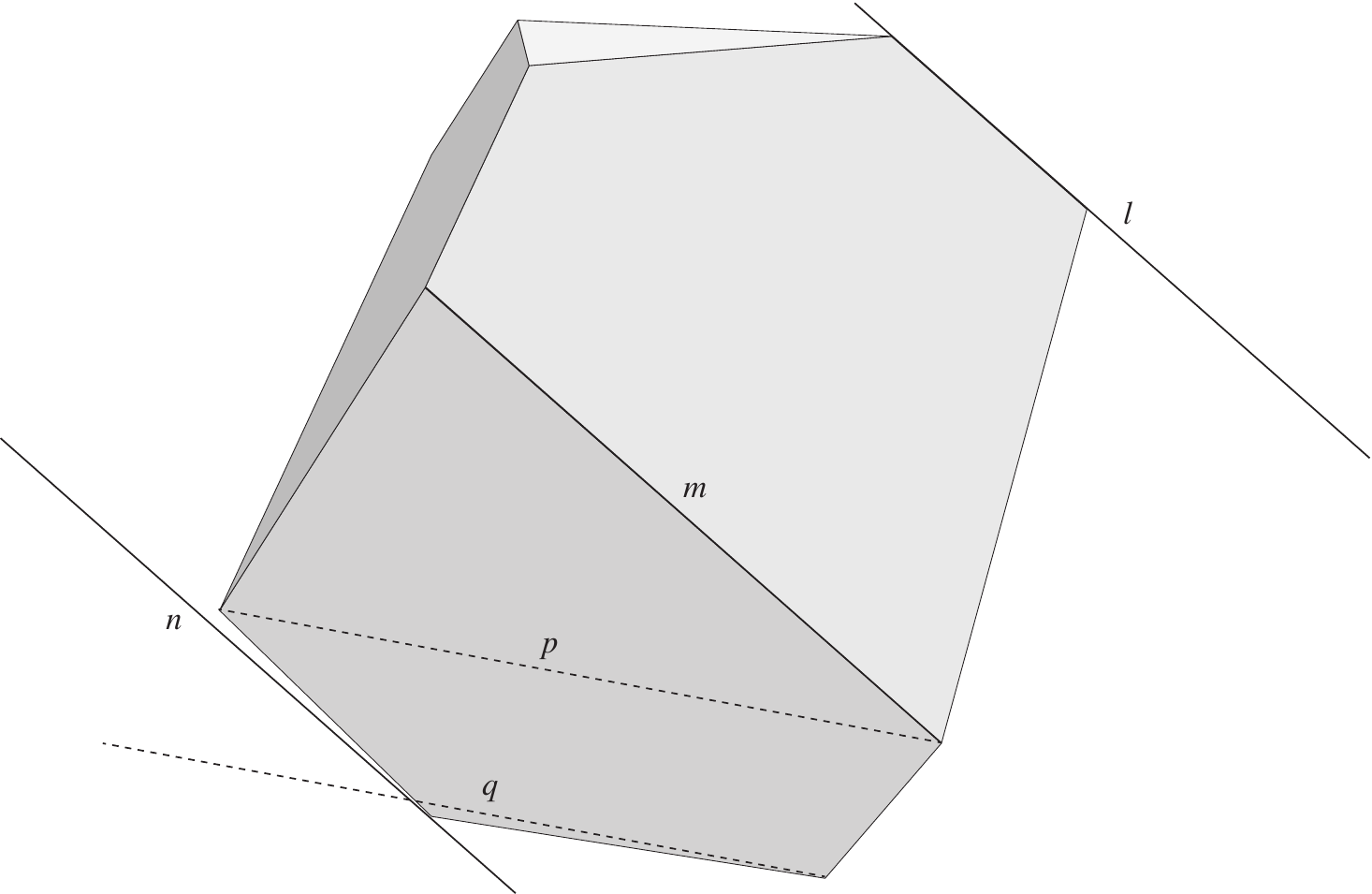}
\caption{Some anomalies in the drawing of the solid.}
\label{anomaly}
\end{figure}

\Du's rendering of the solid has certain anomalous features, a couple of which are depicted in Figure~\ref{anomaly}. The solid line segments $l$, $m$, and $n$ are translations of one another and therefore parallel. Observe that segments $l$ and $m$ very nearly coincide with edges of the solid, while $n$ noticeably diverges from the nearest edge. So let us turn our attention away from the lines in the image to the lines in the (imagined) actual solid object, sitting in space. If the solid is a trigonal trapezohedron as is usually assumed, the three edges that we are considering are parallel in space to one another, hence their images must be concurrent at a vanishing point. But the images nearest $l$ and $m$ are essentially parallel, hence their vanishing point is at infinity (a so-called {\it ideal point\/}). On the other hand, the image of the edge nearest $n$ is clearly not parallel to the images of the other two edges, so the images of the three edge lines are not concurrent at any point, ideal or ordinary. This feature is either an inconsistency in the drawing, or else \Du's conception of the solid contradicts the usual assumptions.

To illustrate another anomaly in Figure~\ref{anomaly}, we have also drawn the dashed line segments $p$ and $q$ parallel to one another. Segment $p$ coincides with the image of a face diagonal, while $q$ is close to but deviates slightly from the image of a bottom edge. Under the usual assumptions about the solid, the diagonal and the edge are parallel in space, hence their images in the engraving should converge to a vanishing point on the horizon line, above and to the left of the picture in Figure~\ref{anomaly}. But notice that as the image of the bottom edge goes to the left, it dips below line $q$, hence it actually diverges from line $p$ (or equivalently, it converges in the wrong direction). The bottom edge should in fact rise above line $q$ (or the diagonal should dip below line $p$) as it goes to the left.

The above errors in the drawing, if they are errors, are not drastic. Indeed, most authors seem not to have noticed them.  However, they highlight the importance of computing the parameters of the truncated rhomboid in multiple ways, using different parts of the image, perhaps even avoiding certain parts if they seem to be too inaccurately drawn.

\end{document}